\documentclass{amsart}

\usepackage{amssymb}

\newtheorem{theorem}{Theorem}[section]
\newtheorem{lemma}[theorem]{Lemma}
\newtheorem{conjecture}[theorem]{Conjecture}
\newtheorem{proposition}[theorem]{Proposition}
\newtheorem{problem}[theorem]{Problem}
\newtheorem{corollary}[theorem]{Corollary}

\theoremstyle{definition}

\newtheorem*{ack}{Acknowledgements}

\theoremstyle{remark}
\newtheorem{remark}[theorem]{Remark}

\numberwithin{equation}{section}

\newcommand \CC {\mathbb C}
\newcommand \ZZ {\mathbb Z}
\newcommand \RR {\mathbb R}

\newcommand \de {\delta}

\newcommand {\cP} {\mathbb {CP}}

\newcommand \C {\mathcal C}
\newcommand \D {\mathcal D}

\newcommand \pa {\partial}

\newcommand \si {\sigma}

\newcommand \cO{\mathcal O}
\newcommand \maxid {\mathfrak m}

\newcommand \Q {\mathbb Q}
\newcommand \spcheck {\vee}
\newcommand \Hom {\textrm{Hom}}
\newcommand \discr {\textrm{discr\,}}
\newcommand \CL {\mathcal L}
\newcommand \rank {\textrm{rank\,}}
\newcommand\CS{{\mathcal S}}
\newcommand\CK{{\mathcal K}}
\newcommand\BA{\bold{A}}
\newcommand\BD{{\bf D}}
\newcommand\BE{{\bf E}}
\newcommand\BU{{\bf U}}
\newcommand\R{{\mathbb R}}
\newcommand\PP{{\mathbb P}}
\newcommand\conj{\textrm{conj}}
\newcommand\Pic{\textrm{Pic}}
\newcommand\supp{\textrm{supp\,}}
\newcommand\Aut{\textrm{Aut}}

\newcommand\CP{\CC\PP}
\newcommand\NN{\mathbb N}
\newcommand\QQ{\mathbb Q}

\begin{document}

\title{On total reality of meromorphic functions}


\author[A.~Degtyarev]{Alex Degtyarev}
\address{Department of Mathematics, Bilkent University,
Bilkent, Ankara 06533, Turkey} \email{degt@fen.bilkent.edu.tr }

\author[T.~Ekedahl]{Torsten Ekedahl}
\address{Department of Mathematics, Stockholm University, SE-106 91
Stockholm, Sweden} \email{teke@math.su.se}

\author[I.~Itenberg]{Ilia Itenberg}
\address{IRMA, Universit\'{e} Louis Pasteur, 7 rue Ren\'{e} Descartes,
67084 Strasbourg Cedex, France} \email{itenberg@math.u-strasbg.fr}
\thanks{I.I. is partially supported by the ANR-network
"Interactions et aspects \'{e}num\'{e}ratifs des
g\'{e}om\'{e}tries r\'{e}elle, tropicale et symplectique"}

\author[B.~Shapiro]{Boris Shapiro}
\address{Department of Mathematics, Stockholm University, SE-106 91
Stockholm, Sweden} \email{shapiro@math.su.se}

\author[M.~Shapiro]{Michael Shapiro}
\address{Department of Mathematics, Michigan State University, East Lansing,
MI 48824-1027, USA} \email{mshapiro@math.msu.edu}
\thanks{M.S. is partially supported by NSF grants DMS-0401178,
by the BSF grant 2002375 and by the Institute of Quantum Science,
MSU}


\date{}

\dedicatory{}

\commby{}

\begin{abstract}
We show that if a meromorphic function of degree at most four on a
real algebraic curve of an arbitrary genus has only real critical
points then it is conjugate to a real meromorphic function  after
a suitable projective automorphism of the image.
\end{abstract}

\maketitle


\section {Introduction}\label{intro}

Let $\gamma:\cP^1\to\cP^n$ be a rational curve in $\cP^n$. We say
that a point $t\in \cP^1$ is a {\em{flattening point}} of $\gamma$
if the osculating frame formed by $\gamma'(t),
\gamma''(t),\dots,\gamma^{(n)}(t)$ is degenerate. In other words,
flattening points of
$\gamma(t)=(\gamma_0(t):\gamma_2(t):\dots:\gamma_{n}(t))$ are
roots of the Wronskian $$W(\gamma_0,\dots,\gamma_{n})=
\left|%
\begin{array}{ccc}
  \gamma_0 & \dots & \gamma_{n} \\
  \gamma_0' & \dots & \gamma_{n}' \\
  & \dots & \\
  \gamma_0^{(n)} & \dots & \gamma_{n}^{(n)} \\
\end{array}%
\right|. $$

In 1993  B.~ and M.~Shapiro made the following claim which we will refer to as
{\em rational total reality conjecture}.
\begin{conjecture}\label{conj1}
 If all flattening
points of a rational curve $\cP^1\to\cP^n$ lie on the real line
$\RR P^1\subset \cP^1$ then the curve is conjugate to a real
algebraic curve under an appropriate  projective automorphism  of
$\cP^n$.
\end{conjecture}
Notice that coordinates $\gamma_i$ of the rational curve $\gamma$ are homogeneous
polynomials of a certain degree, say $d$. Considering them as vectors in the
space of homogeneous degree $d$ polynomials   we can reformulate
the above conjecture as a statement  of total reality in Schubert
calculus, see  ~\cite{KhS},~\cite{So1}-\cite{So6},~\cite{Ve}. Namely, for any $0\le d<n$
let $t_1<t_2<\dots<t_{(n+1)(d-n)}$ be a sequence of real numbers and
$\bf {r}: \CC\to \CC^{d+1}$ be a rational normal curve with coordinates
 $r_i(t)=t^i,\, i=\overline{0, d}$. Denote by  $T_i$  the
osculating $(d-n)$-dimensional plane to $\bf{r}$ at the moment $t=t_i$.
Then the above rational total reality conjecture is equivalent to the following claim.
\begin{conjecture} [Schubert calculus interpretation]\label{conj2}
In the above notation  any $(n+1)$-dimensional subspace in $\CC^{d+1}$
 which  meets all $(n+1)(d-n)$ subspaces $T_i$ nontrivially  is real.
\end{conjecture}
It was first supported by extensive numerical evidences, see \cite
{So1}-\cite{So6},~\cite {Ve} and later settled for $n=1$, see
\cite{EG}. The case $n\ge 2$ resisted all efforts for a long time.
In fall 2005 the authors were informed by   A.~Eremenko and
A.~Gabrielov that they were able to prove  Conjecture~\ref{conj1} for plane rational quintics.
Just few months later it  was completely established by
E. Mukhin, V. Tarasov, and A. Varchenko
in~\cite{MukhinTarasovVarchenko}.

Their  proof reveals the deep
connection between Schubert calculus and theory of integrable
system and is based on the Bethe ansatz method in the Gaudin
model. More exactly, conjectures 1 and 2  are reduced to the question of reality
of $(n+1)$-dimensional subspaces of the space $V$ of polynomials of
degree $d$ with given asymptotics at infinity and fixed Wronskian.
Choosing a base  in such a subspace we get the rational curve
 $\cP^1\to \cP^n$, whose flattening points coincide with
the roots of the above mentioned Wronskian. The subspaces with
desired properties are constructed explicitly using properties of
spectra of Gaudin Hamiltonians. Namely, relaxing the reality
condition these polynomial subspaces are labeled by common
eigenvectors of Gaudin Hamiltonians, one-parameter families of
commuting linear maps on some vector space,
   $H_1(x),\dots,H_{n+1}(x) : V\to V$.
The subspace, labeled by an  eigenvector, is the kernel of a
certain linear differentail operator of order $n+1$, assigned to
each eigenvector of the Hamiltonians. The  coefficents of that
differentail operator, are the eigenvalues of the Hamiltonians on
that eigenvector.
%
It turns out that in the
case of real rooted
 Wronskians  Gaudin Hamiltonians are
symmetric with respect to the so-called tensor Shapovalov form, and  thus have real spectra. Moreover,
their eigenvalues are real rational functions. This fact implies that the  kernels
of the above fundamental differential operators are real subspaces in $V$ which
concludes the proof.

Meanwhile two different generalizations   of the original conjectures  (both dealing with  the case $n=1$)  were suggested in~\cite{EGSV}
and~\cite{ESS}. The former replaces
the condition of reality of critical points by the existence of
separated collections of real points such that a meromorphic
function takes the same value on each set.
 The latter discusses the generalization
of the total reality conjecture to higher genus curves.

The present paper is the sequel of~\cite{ESS}. Here we prove the higher
genus version of the total reality conjecture for all meromorphic
funtions of degree at most four.



For reader's convenience and to make the paper self-contained we
included some of results of~\cite{ESS} here. We start with some standard notation.


\medskip
\noindent {\bf Definition.} A pair $(\C,\sigma)$ consisting of a
compact Riemann surface $\C$ and its antiholomorphic involution
$\sigma$  is called a {\em real algebraic} curve. The set
$\C_{\sigma}\subset \C$ of all fixed points of $\sigma$ is called
the {\em real part} of $(\C,\sigma)$.

If $(\C,\sigma)$ and
$(\D,\tau)$ are real curves (varieties) and ${f}:{\C}\to{\D}$
a holomorphic map, then we denote by  $\overline{f}$ the holomorphic map $\tau\circ
f\circ\sigma$. Notice  that  $f$ is {\em real}
if  and only if  $\overline{f}=f$.

The main question we discuss  below is as follows.

\medskip
\noindent {\bf Main Problem.} Given a meromorphic function
$f:(\C,\sigma)\to \cP^1$ such that

\noindent
i) all its critical points and values are distinct;

\noindent
 ii) all its critical points belong to $\C_{\sigma}$;

\noindent is it true that that $f$ becomes a real meromorphic
function after an appropriate  choice of a real structure on $\cP^1$?


\medskip
\noindent {\bf Definition.} We say that the space of meromorphic
functions of degree $d$ on a genus $g$ real algebraic curve
$(\C,\sigma)$ has \emph{the total reality property} (or is
{\emph{totally real}) if the Main Problem has the affirmative
answer for any meromorphic function from this space which
satisfies the above assumptions. We say that a pair of positive integers $(g,d)$
has a total reality property if the space of meromorphic functions of degree $d$ is totally real on
{\bf any} real algebraic curve of genus $g$.

\medskip
The following results were proven in~\cite{ESS} (see Theorem~1 and
Corollary~1 there).

 \begin{theorem}  The space of  meromorphic  functions of  any degree $d$
 which is a prime on any real curve $(\C,\sigma)$ of genus $g$ which
 additionally satisfies the inequality:
 $g>\frac{d^2-4d+3}{3}$  has the total reality property.
\label{th:sqrt}
\end{theorem}

\begin{corollary}  The  total reality property holds for all meromorphic functions
 of  degrees $2,3$, i.e. for all pairs $(g,2)$  and $(g,3)$.  \label{th:main}
\end{corollary}

The proof of the Theorem~\ref{th:sqrt} is based on the following
observation. Consider the space $\cP^1\times \cP^1$ equipped with the involution
$s: (x,y)\mapsto(\bar y,\bar x)$ which we call  the {\em involutive
real structure} (here $\bar x$ and $\bar y$ stand for  the complex
conjugates of $x$ and $y$  with respect to the standard real structure in
$\cP^1$). The pair $ \bf{Ell} = (\cP^1\times \cP^1,s)$   is usually referred
 to as the standard  {\em ellipsoid}, see \cite{GSh}.  (Sometimes by the ellipsoid one means
 the set of fixed points of $s$ on $\cP^1\times \cP^1$). The next statement translates the
 problem of total reality into the question of  (non)existence of certain real
 algebraic curves on $\bf{Ell}$.

\begin{proposition}\label{prop:Equiv.Curve} For any positive integer $g$ and {\bf prime} $d$
the total reality property holds for the pair $(g,d)$  if and only if there is
{\bf no} real algebraic curve on $\bf{Ell}$ with the following properties:

\noindent
{\em i)}   its geometric genus equals $g$;

\noindent
{\em ii)} its  bi-degree as a curve on $\cP^1\times \cP^1$ equals $(d,d)$;

\noindent
{\em  iii)} its only singularities are   $2d-2+2g$ real cusps
on $\bf {Ell}$ and possibly some number of (not necessarily transversal)
intersections of smooth branches.
\end{proposition}

Extending slightly the arguments proving Proposition~\ref{prop:Equiv.Curve} one gets the following
statement.

\begin{proposition}\label{prop:General}  The total reality property  holds for all real meromorphic functions,
i.e. for all pairs $(g',d')$ if and only if for {\bf no} pair
$(g,d),\ d>1$ there exists a real algebraic curve on   $\bf{Ell}$
satisfying conditions  {\em i)} - {\em iii)} of
Proposition~\ref{prop:Equiv.Curve}.
\end{proposition}

The main result of the present paper obtained using a version of Proposition~\ref{prop:Equiv.Curve} and technique related to integer lattices and $K3$-surfaces is as follows.

 \begin{theorem}\label{thm:2}
 The  total reality property holds for all meromorphic functions
 of  degree $4$, i.e. for all pairs $(g,4)$.
\end{theorem}


\medskip
The structure of the note is as follows. \S~\ref{sc:pr} contains
the proofs of Theorem~\ref{th:sqrt}, Corollary~\ref{th:main} and
reduction of Theorem~\ref{thm:2} to the question of nonexistence
of a real curve $\D$  on $\bf{Ell}$ of bi-degree $(4,4)$ with eight real cusps
and no other singularities.
\S~\ref{sc:Itenberg} contains a proof of nonexistence of such
curve $\D$,  while \S~\ref{sc:rmk} contains a number of remarks
and open problems.

\begin{ack}
The authors are  grateful to A.~Gabrielov, A.~Eremenko,
R.~Kulkarni, B.~Osserman, V.~Tarasov, A.~Vainshtein, and
A.~Varchenko
 for  discussions of the topic.
The third, fourth and fifth authors want to
acknowledge the hospitality of MSRI in Spring 2004 during the
program 'Topological methods in real algebraic geometry' which
gave them a large number of valuable research inputs.
    \end{ack}

\section {Proofs}
\label{sc:pr}

 If not mentioned explicitly we assume below  that $\cP^1$ is
provided with its standard real structure.  Real meromorphic functions on a real algebraic curve $(\C,\sigma)$
can be characterized in the following way.

\begin{proposition}\label{realification criterion}
If $(\C,\sigma)$ is a proper irreducible real curve and
$f:{\C}\to{\cP^1}$ is a non-constant holomorphic map, then $f$ is real for
some real structure on $\cP^1$ if and only if  there is a M\"{o}bius
transformation ${\varphi}:{\cP^1}\to{\cP^1}$ such that
$\overline{f}=\varphi\circ f$.

\begin{proof}
Any real structure on $\cP^1$ is of the form $\tau\circ\phi$ for a
complex M\"{o}bius transformation $\phi$ and $\tau$ the standard
real structure with $\overline{\phi}\circ\phi=\text{id}$ and
conversely any such $\phi$ gives a real structure. If $f$ is real
for such a structure we have $f=\tau\circ \phi\circ f\circ\sigma$,
i.e., $\overline{f}=\varphi\circ f$ for
$\varphi=\tau\circ\phi^{-1}\circ\tau$. Conversely, if
$\overline{f}=\varphi\circ f$, then $f=
\overline{\overline{f}}=\overline{\varphi\circ
f}=\overline{\varphi}\circ\overline{f}=\overline{\varphi}\circ\varphi\circ
f$ and as $f$ is surjective we get
$\overline{\varphi}\circ\varphi=\text{id}$. That means that
$\phi:=\tau\circ\varphi^{-1}\circ\tau$, then $\phi$ defines a real
structure on $\cP^1$ and by construction $f$ is real for that
structure and the fixed one on $\C$.
\end{proof}
\end{proposition}

Up to  a real isomorphism there are only two real structures on
$\cP^1$, the standard one and the one on an isotropic real quadric
in $\cP^2$. The latter is distinguished from the former by not
having any real points.

Assume now that $(\C,\sigma)$ is a proper irreducible real curve and ${f}:{\C}\to{\cP^1}$ a
non-constant  meromorphic function. It defines the holomorphic map
\begin{displaymath}
\C \stackrel{(f,\overline{f})}{\longrightarrow} \cP^1\times\cP^1
\end{displaymath}
and if $\cP^1\times\cP^1$ is given the involutive real structure
 $s: (x,y) \to (\bar y,\bar x)$
then it is clearly a real map. Now we can formulate the central technical result of this section.

\noindent
\begin{proposition}\label{real diagonal} \

\begin{enumerate}
\item \label{it1} The image $\D$ of the curve $\C$ under the map $(f,\overline{f})$ is of type
$(\de,\de)$ for some positive integer $\de$ and if $\pa$ is the degree of the map $\C \to
\D$ we have that $d=\de\pa$, where $d$ is the degree of the original $f$.

\item \label{it2} The function $f$ is real for some real structure
on $\cP^1$ precisely when $\de=1$.

\item \label{it3} Assume that $\C$ is smooth and all the critical points of $f$ are real. Then all the
critical points of  ${\psi}:{\widetilde \D}\to{\cP^1}$, the composite of the normalization
map $\widetilde \D \to \D$ and the restriction of the projection of $\cP^1\times\cP^1$
has all its critical points real.
\end{enumerate}
\end{proposition}

\begin{proof}
The image of $\C$  under the real holomorphic map $(f,\overline{f})$  is a real curve so that
$\D$ is a real curve in $\cP^1\times\cP^1$ with respect to  its involutive real
structure, i.e. a real curve on the ellipsoid $\bf{Ell}$. Any such curve is of type $(\de,\de)$ for some positive integer  $\de$  since  the involutive  real structure permutes the
two degrees. The rest of (\ref{it1}) follows by using the multiplicativity of degrees for the
maps $f\colon \C \to \D \to \cP^1$, where the last map is projection on the first
factor.

As for (\ref{it2}) assume first that $f$ can be made real for some real
structure on $\cP^1$. By Proposition \ref{realification criterion} there is a
M\"{o}bius transformation $\varphi$ such that $\overline{f}=\varphi\circ f$ but that
in turn means that $(f,\overline{f})$ maps $\C$ into the graph of $\varphi$ in
$\cP^1\times\cP^1$.  This  graph is hence equal to $\D$ and is thus of type
$(1,1)$. Conversely, assume that $\D$ is of type $(1,1)$. Then it is a graph of
an isomorphism $\varphi$ from $\cP^1$ to $\cP^1$ and by construction
$\overline{f}=\varphi\circ f$ so we conclude by another application of
Proposition \ref{realification criterion}.

Finally, for (\ref{it3}) we have that the map $\C \to \D$ factors as a,
necessarily real, map ${h}:{\C}\to{\widetilde \D}$ and then $f=\psi \circ h$. If $pt \in
\widetilde \D$ is a critical point, then all points of $h^{-1}(pt)$ are critical for
$f$ and hence by assumption real. As $h$ is real this implies that $pt$ is also
real.
\end{proof}

Part (\ref{it2}) of the above Proposition gives another reformulation
of the  total reality property  for meromorphic functions.

\begin{corollary}\label{cor:reform}  If a  meromorphic  function $f:(\C,\si)\to \cP^1$ of degree $d$
is real for some real structure on $\cP^1$ then the map $\C
\stackrel{(f,\overline{f})}{\longrightarrow} \D\subset \cP^1\times
\cP^1$ must have degree $d$ as well.
\end{corollary}

\begin{remark} Notice that without the requirement of reality of $f$ the degree of
$\C \stackrel{(f,\overline{f})}{\longrightarrow}  \D$ can be any factor of $d$.
\end{remark}

By a {\em cusp} we  mean a  curve singularity of multiplicity $2$
and whose tangent cone is a double line. It has the local form
$y^2=x^{k}$ for some integer $k\ge 3$ where $k$ is an invariant
which we shall call its {\em type}. A cusp of type $k$ gives a
contribution of $\lceil(k-1)/2\rceil$ to the arithmetic genus of a
curve. A cusp of type $3$ will be called {\em ordinary}.

If $\C$ is a curve and $p_1,\dots,p_k$ are its smooth points then
consider the finite map  ${\pi}:{\C}\to{\C(p_1,\dots,p_k)}$ which
is a homeomorphism and for which $\cO_{\C(p_1,\dots,p_k)} \to
\pi_*\cO_\C$ is an isomorphism outside of $\{p_1,\dots,p_k\}$ such
that the image of the map
$\cO_{\C(p_1,\dots,p_k),\pi(p_i)}\to\cO_{\C,p_i}$ is the inverse
image of $\CC$ in $\cO_{\C,p_i}/\maxid_{p_i}^2$. In other words,
$\C(p_1,\dots,p_k)$ has ordinary cusps at all points $\pi(p_i)$.

Then $\pi$ has the following two (obvious) properties:


\begin{lemma}\label{lm:properties} \

\begin{enumerate}
\item A holomorphic map ${f}:{\C}\to{X}$ which is not an immersion at all the points
$p_1,\dots,p_k$ factors through $\pi$.

\item If $\C$ is proper, then the arithmetic genus of $\C(p_1,\dots,p_k)$ is $k$
plus the arithmetic genus of $\C$.
\end{enumerate}
\end{lemma}

\begin{proposition}\label{real classification}
Assume that $(\C,\si)$ is a smooth and proper real curve and let
${f}:~{\C}\to{\cP^1}$ be a meromorphic function of degree $d$. If there
are $k$ real points $p_1,\dots,p_k$ on $\C$ which are critical
points for $f$ and if $(f,\overline{f})$ gives a map of degree $1$
from $\C$ to its image $\D$ in $\cP^1\times\cP^1$, then
$g(\C)+k\le (d-1)^2$. If $g(\C)+k=(d-1)^2$, then the map $h:\C \to
\D$ factors to give an isomorphism $\C(p_1,\dots,p_k)
\stackrel{\sim}{\rightarrow} D$.
\begin{proof}
As $p_i$ is a real critical point  it is a critical  for $\overline{f}$ as well  and hence for
$(f,\overline{f})$. This implies by the first property for $\pi: \C \to
\C(p_1,\dots,p_k)$ that the map $\C \to D$ factors as $\C \to \C(p_1,\dots,p_k) \to
\D$ and hence the arithmetic genus of $\C(p_1,\dots,p_k)$, which is $g(\C)+k$ by
the second property of $\C(p_1,\dots,p_k)$, is less than or equal to the
arithmetic genus of $\D$, which by the adjunction formula is equal to
$(d-1)^2$. If we have equality then their genera are equal and hence the map
$\C(p_1,\dots,p_k) \to \D$ is an isomorphism.
\end{proof}
\end{proposition}

Now it is easy to derive Proposition~\ref{prop:Equiv.Curve} from
Proposition~\ref{real diagonal}}. Indeed,  if  a meromorphic function $f:\C\to \cP^1$ of a prime degree $d$
with all $2g+2d-2$ real critical points
  can not be made real then its image under
$(f,\bar f)$ in $\cP^1\times\cP^1$   is the real curve on $\bf{Ell}$  with
$2g+2d-2$ real cusps and no other  singularities different from intersections of smooth branches.   (Intersections of smooth branches in the image might occur and are moreover necessary to produce the required genus.)  Vice versa,
assume that such a  curve $\D\subset \cP^1\times\cP^1$ which is real in the involutive structure does exist.
Let $\widetilde{\D}$ be the normalization of $\D$, and consider the natural birational  projection map $\mu:\widetilde{\D}\to
\D$. Define $f:\widetilde{\D}\to\cP^1$ as a
composition $\widetilde{D}\to \D\to\cP^1$, where the last map is induced by
the projection of $\cP^1\times\cP^1$ on the first factor. It
remains to notice that all $2g+2d-2$ critical points of $f$ are
real while $f$ can not be made real by Proposition~\ref{real
diagonal}.\qed

Similar arguments show the validity of
Proposition~\ref{prop:General}. Indeed, assume  that  there exists a
meromorphic function $\phi$ of some degree $d'$ on a real curve
$\C'$ of some genus $g'$ violating the total reality conjecture.  Let   $\D' \subset \cP^1\times \cP^1$
be  the image curve  of bidegree
$(d,d)$ obtained by application of the map $(\phi,\bar \phi)$ to
$\C'$ and  let $\widetilde{\D'}$ be the normalization of $\D'$. Let
 $\mu': \widetilde{\D'} \to\D'$ be the canonical birational map and, finally, let $\phi:\widetilde{\D'}\to\CP^1$
be the composition of $\mu'$ and the projection of $\cP^1\times\cP^1$
on its first factor. Then $\phi$ has degree $d$ and all the
critical points of $\phi$ are real by
Proposition~\ref{real diagonal}(3). Note that if $f$ is not
conjugate to a real function by a M\"obius transformation the same holds for $\phi$ as well.
Hence, $\phi:\widetilde{D'}\to\cP^1$ also violates the total
reality conjecture. The image of
$\widetilde{D'}\stackrel{(\phi,\overline{\phi})}{\longrightarrow}\cP^1\times\cP^1$
coincides with $\D'$, and  the map $\mu':\widetilde{D'}\to\D'$ is birational. So $\D'$
satisfies the assumptions i)-iii) of
Proposition~\ref{prop:Equiv.Curve} for $g=g(\widetilde{D'})$ and $d=\delta>1$,
see Proposition~\ref{real diagonal}. Indeed, the map $\C'\to\D'$
lifts to a map $\C'\to\widetilde{D'}$ of degree $\delta=d'/d$ with only simple
ramifications whose number by the Riemann-Hurwitz formulas is
$2g(\C')-2-\delta(2g(\widetilde{D'})-2)$. Hence the number of critical points
of $f$ that are the preimages of cusps of $\D'$ can be computed as
$K=2g(\C')-2+2d'-(2g(\C')-2-\delta(2g(\widetilde{D'})-2))$. Note that each cusp
has as preimages exactly $\delta$ critical points. Finally we
compute the number of cusps of $\D'$ as
$\frac{1}{\delta}K=2g(\widetilde{D'})-2+2d$.

And conversely, exactly as in the above proof given a curve
$\D' \subset \cP^1\times \cP^1$ satisfying the assumptions i)-iii)
of Proposition~\ref{prop:Equiv.Curve}  we get a meromorphic
function violating the total reality conjecture by composing the
birational projection $\mu'$ from the normalization $\widetilde
{\D'}$  to $\D'$ with the  projection of $\D'$ on the
first coordinate in $\cP^1\times \cP^1$. \qed

Now we are ready to prove Theorem~\ref{th:sqrt}.
 Indeed,
 if we assume that all the critical points of a generic meromorphic
 function $f:(\C,\sigma)\to\cP^1$
are real then $k$ in the above Proposition equals $2d-2+2g(\C)$.
Under the assumption $g(\C)>\frac{d^2-4d+3}{3}$ one gets
$g(\C)+k=2d-2+3g(\C)>(d-1)^2$. Thus, the case when $\C$ maps
birationally to $\D$ is impossible by Proposition~\ref{real classification}. Since
$d$ is prime the only other possible case is when the degree of
the map $\C\to \D$ equals $d$ and therefore, the degree of the map
$\D\to \cP^1$ equals $1$ which by (\ref{it2}) of
Proposition~\ref{real diagonal} gives that the total reality
property holds. \qed

\medskip
Let us  apply Theorem~\ref{th:sqrt} to prove Corollary~\ref{th:main}.

\medskip
\noindent Case $d=2$. Suppose that the degree $d$ of a meromorphic function
${f}:(\C,\si)\to{\cP^1}$ is equal to $2$. That only leaves two
possibilities: The first is that the map $\C \to \D$ has degree
$2$ and then by Proposition \ref{real diagonal} $f$ is real for
some real structure $(\cP^1,\tau)$ on $\cP^1$. In particular, if
the set $\C_\si$ of real points is nonempty then $(\cP^1,\tau)$
has the same property which means that it is equivalent to the
standard real structure. The second is that the map $\C \to \D$ is
birational and then by Proposition \ref{real classification} we
get $g(\C)+k \le 1^2=1$, where $k$ is the number of real critical
points of $f$. In particular if $g(\C)>0$ then there are no real
critical points. Thus a hyper-elliptic map from a real curve
$(\C,\si)$ is real if one of its critical points is real.

\noindent
Case $d=3$.  In this case again we have only two possibilities; either $f$
is real for a real structure on $\cP^1$ or $\C \to \D$ is birational in which case
we have $g(\C)+k \le 2^2=4$. The case $g(\C)=0$ was settled in \cite {EG}. Recall that
the total number of critical points equals $2d-2+2g(\C)$. But if $g(\C)>0$
then $2\cdot 3-2+3g>4$ and this case of Theorem~\ref{th:main} is settled.
Analogously to the case $d=2$  a function $f$  with the degree  $d=3$ is real if it has more
than $\max(4-g(\C),1)$ real critical points. \qed

\medskip
Now we can finally start proving Theorem~\ref{thm:2}. Using a
version of Proposition~\ref{real diagonal} we reduce the case of
degree $d=4$  to the existence problem of a real curve on the
ellipsoid $\bf{Ell}= (\cP^1\times\cP^1,s)$ of bi-degree $(4,4)$
with 8 ordinary real cusps and {\bf no} other singularities.
Indeed, arguing along the same lines as above  we have three
possibilities for the image $\D$ of $\C$ under the map $(f,\bar
f)$. Namely, $\D$ might have bi-degrees $(1,1)$, $(2,2)$, or
$(4,4)$. In the first case $f$ can be made real. In the second
case, by Proposition \ref{real diagonal}, the projection on the
first factor will give a map from the normalization $\widetilde
\D$ of $\D$. The arithmetic genus $p_a(\D)=1$, and the geometric
genus $g(\widetilde \D)$ of the normalization $\widetilde D$ does
not exceed $1$. Let $\widetilde h:\C\to\widetilde\D$ be the lift
of $h:\C\to\D$. Note that if $p_i\in \C$ is a critical point of
$f$ then either its image $h(p_i)$ is a cusp of $\D$ or $p_i$ is a
ramification point of $\widetilde h$. The ramification divisor
$R(\widetilde h)=2g(\C)+2-4g(\widetilde\D)$. The number of cusps
of $\D$ does not exceed $1$, whereas the number of distinct
critical points of $f$ is $2g(\C)+6$. Note that any cusp has two
critical points of $f$ as preimages. Therefore, we must have
$\frac{1}{2}\left(2g(\C)+6-\left(2g(\C)+2-4g(\widetilde\D)\right)\right)\le
1$ which is impossible.

 We are hence left with the case when $\D$ has degree
$(4,4)$. The only case when $2\cdot 4-2+3g(\C)\le 9$ for $g(\C)>0$
is the case  of $g(\C)=1$. If all the critical points
$p_1,\dots,p_8$ of  ${f}:{\C}\to{\cP^1}$ are real, then we get a
birational map $\C(p_1,\dots,p_8) \to \D$ and as then both
$\C(p_1,\dots,p_8)$ and $\D$ have arithmetic genus $9$, this map
is an isomorphism. Hence $\D$ is a curve with $8$ ordinary real
cusps and no other singularities.  To finish the proof of Theorem~\ref{thm:2} we have to show that such curves do not exist.   Since the proof of this claim  requires a lot of additional notation and techniques we
decided to place it into a separate section.


\section{Intersection lattices}\label{sc:Itenberg}

We will need a number of standard notions from the lattice theory and $K3$-surfaces.

\subsection{Discriminant forms}\label{par:0}

A {\it lattice\/} is a finitely generated free abelian group~$L$
supplied with a symmetric bilinear form $b:L\otimes L\to\ZZ$. We
abbreviate $b(x,y)=x\cdot y$ and $b(x,x)=x^2$. A lattice~$L$ is
{\it even\/} if $x^2=0\bmod2$ for all $x\in L$.
As the transition matrix between two integral bases
has determinant $\pm1$, the
determinant $\det L\in\ZZ$
({\it i.e.}, the determinant
of the Gram matrix of~$b$ in
any
basis of~$L$)
is well defined.
A lattice~$L$ is called
{\it nondegenerate\/} if the determinant $\det L\ne0$; it is called {\it
unimodular\/} if $\det L=\pm1$.


Given a lattice~$L$,
the bilinear form can be extended to $L\otimes\Q$ by linearity. If
$L$ is nondegenerate, the dual group $L^\spcheck=\Hom(L,\ZZ)$ can
be identified with the subgroup
$$
\bigl\{x\in L\otimes\Q\bigm|
 \text{$x\cdot y\in\ZZ$ for all $x\in L$}\bigr\}.
$$
In particular, $L\subset L^\spcheck$. The quotient $L^\spcheck/L$
is a finite group; it is called the {\it discriminant group\/}
of~$L$ and is denoted by $\discr L$ or~$\CL$. The discriminant
group~$\CL$ inherits from $L\otimes\Q$ a symmetric bilinear form
$\CL\otimes\CL\to\Q/\ZZ$,
called the {\it discriminant form},
and, if $L$ is even, its quadratic
extension $\CL\to\Q/2\ZZ$.
When
speaking about the discriminant groups, their
(anti-)isomorphisms, etc., we always assume that the discriminant
form (and its quadratic extension if the lattice is even) is taken
into account. One has $\#\CL=\mathopen|\det L\mathclose|$; in
particular, $\CL=0$ if and only if $L$ is unimodular.


In what follows we denote by~$\BU$ the {\it hyperbolic plane},
{\it i.e.},
the lattice generated by a pair of vectors~$u$,~$v$
(referred to as a {\it standard basis\/} for~$\BU$)
with
$u^2=v^2=0$ and $u\cdot v=1$. Furthermore, given a lattice~$L$, we
denote by~$nL$, $n\in\NN$, the orthogonal sum of $n$~copies
of~$L$, and by~$L(p)$, $p\in\QQ$, the lattice obtained from~$L$ by
multiplying the form by~$q$ (assuming that the result is still an
integral lattice). The notation $n\CL$ is also used for the
orthogonal sum of $n$~copies of a discriminant group~$\CL$.

Two lattices~$L_1$, $L_2$ are said to have the same {\it genus} if
all localizations $L_i\otimes\QQ_p$, $p$~prime, and $L_i\otimes\Q$
are pairwise isomorphic. As a general rule, it is relatively easy
to compare the genera of two lattices; for example, the genus of
an even lattice is determined by its signature and the isomorphism
class of the discriminant
group,
see~\cite{Nikulin}. In the same
paper~\cite{Nikulin} one can find a few classes of lattices whose
genus is known to contain a single isomorphism class.

Following V.~V.~Nikulin, we denote by $\ell(\CL)$ the minimal
number of generators of a finite group~$\CL$ and, for a prime~$p$,
let $\ell_p(\CL)=\ell(\CL\otimes\ZZ_p)$. (Here $\ZZ_p$ stands for
the cyclic group $\ZZ/p\ZZ$.) If $L$ is a nondegenerate lattice,
there is a canonical epimorphism
$\Hom(L,\ZZ_p)\to\CL\otimes\ZZ_p$. It is an isomorphism if and
only if $\rank L=\ell_p(\CL)$.

An {\it extension\/} of a lattice~$L$ is another lattice~$M$
containing~$L$. An extension is called {\it primitive\/} if
$M/L$
is torsion free. In what follows we are only interested in the
case when both~$L$ and~$M$
are even. The relation between
extensions of even lattices and there discriminant forms was
studied in details by Nikulin; next two theorems are found
in~\cite{Nikulin}.

\begin{theorem}\label{thm:Nik1} Given a nondegenerate even
lattice~$L$, there is a canonical one-to-one correspondence
between the set of isomorphism classes of finite index extensions
$M\supset L$ and the set of isotropic subgroups $\CK\subset\CL$.
Under this correspondence one has
$M=\bigl\{x\in
L^\spcheck\bigm|x\bmod L\in\CK\bigr\}$ and $\discr
M=\CK^\perp/\CK$.
\end{theorem}

\begin{theorem}\label{thm:Nik2} Let
$M\supset L$ be a primitive
extension of a nondegenerate even lattice~$L$ to a unimodular even
lattice~$M$.
Then there is a canonical anti-isometry $\CL\to\discr
L^\perp$ of discriminant forms; its graph is the kernel
$\CK\subset\CL\oplus\discr L^\perp$ of the finite index extension
$M\supset L\oplus L^\perp$, see Theorem~\ref{thm:Nik1}.
Furthermore, a pair of auto-isometries of~$L$ and~$L^\perp$
extends to an auto-isometry of~$M$
if and only if the induced
automorphisms of~$\CL$ and $\discr L^\perp$, respectively,
agree via the above
anti-isometry of the discriminant
groups.
\end{theorem}

The general case
$M\supset L$ splits into the finite index
extension $\tilde L\supset L$ and primitive extension
$M\supset\tilde L$, where
$$
\tilde L=\bigl\{x\in
M\bigm|nx\in L\ \text{for some $n\in\ZZ$}\bigr\}
$$
is the {\it primitive hull\/} of~$L$ in~$M$.

A {\it root\/} in an even lattice~$L$ is a vector $r\in L$ of
square~$-2$. A {\it root system\/} is an even negative definite
lattice generated by its roots. Recall that each root system
splits (uniquely up to order of the summands) into orthogonal sum
of indecomposable root systems, the latter being those of types
$\BA_p$, $p\ge1$, $\BD_q$, $q\ge4$, $\BE_6$, $\BE_7$, or~$\BE_8$,
see~\cite{Bourbaki}. A finite index extension
$\Sigma\subset\tilde\Sigma$ of a root system~$\Sigma$ is called
{\it quasi-primitive\/} if each root of~$\tilde\Sigma$ belongs
to~$\Sigma$.

Each root system that can be embedded in~$\BE_8$ is unique in its
genus,
see~\cite{Nikulin}.
In what follows we need the discriminant group
$\discr\BA_2=\left<-\frac{2}{3}\right>$: it is the cyclic group
$\ZZ_3$ generated by an element of square $-\frac{2}{3}\bmod2\ZZ$.

\subsection{$K3$-surfaces and ramified double coverings of
$\CP^1 \times \CP^1$}\label{par:1}
A {\it $K3$-surface} is a nonsingular compact connected and simply
connected complex surface with trivial first Chern class.
From the Castelnuovo--Enriques classification of surfaces it
follows that all $K3$-surfaces form a single deformation family.
In particular, they are all diffeomorphic, and the calculation for
an example
(say, a quartic in~$\CP^3$)
shows that
$$
\chi(X)=24,\quad h^{2,0}(X)=1,\quad h^{1,1}(X)=20.
$$
(see, for instance, \cite{BPV}).
Hence, the intersection lattice $H_2(X;\ZZ)$ is an even (since
$w_2(X)=K_X\bmod2=0$) unimodular (as intersection lattice of any
closed $4$-manifold) lattice of rank~$22$ and signature~$-16$. All
such lattices are isomorphic to $L=2\BE_8\oplus3\BU$.
In
particular, the quadratic space $H_2(X;\R)\cong L\otimes\R$ has
three positive squares; for a maximal positive definite subspace
one can choose the subspace spanned by the real and imaginary
parts of the class~$[\omega]$ of a holomorphic form~$\omega$
on~$X$ and the class~$[\rho]$ of the fundamental form of a
K\"{a}hler metric on~$X$. (We identify the homology and cohomology
{\it via\/} the Poincar\'{e} duality.)

A {\it real} $K3$-surface is a pair $(X, \conj)$, where~$X$ is a
$K3$-surface and $\conj: X \to X$ an anti-holomorphic involution.,
i.e., a real structure on~$X$. The $(+1)$-eigenlattice $\ker(1 -
\conj_*) \subset H_2(X; \ZZ)$ of~$\conj_*$ is hyperbolic, i.e., it
has one positive square in the diagonal form over~$\R$. This
follows, e.g., from the fact that $\omega$ and~$\rho$ above can be
chosen so that $\conj_*[\omega]=[\bar\omega]$ and
$\conj_*[\rho]=-[\rho]$.

Let $Y=\CP^1\times\CP^1$ and let $C\subset Y$ be an
irreducible curve of bi-degree $(4, 4)$ with at worst simple
singularities ({\it i.e.}, those of type~$\BA_p$, $\BD_q$,
$\BE_6$, $\BE_7$, or~$\BE_8$). Then, the minimal resolution~$X$ of
the double covering of~$Y$ ramified along~$C$ is a $K3$-surface.
Recall that the standard ellipsoid is the pair ${\bf Ell}=(Y,s)$ where
$s$ is the anti-holomorphic involution
$s:Y\to Y$,
$s:(x,y)\mapsto (\bar y,\bar x)$.
If~$C$ is $s$-invariant, the involution~$s$ lifts to two different
real structures on~$X$, which commute with each other and with the
deck translation of the covering $X\to Y$. Choose one of the two
lifts and denote it by $\conj$.

Let $l_1,l_2\in H_2(X;\ZZ)$ be the pull-backs of the classes of
two lines belonging to the two rulings of~$Y$. Then
$l_1^2=l_2^2=0$ and $l_1\cdot l_2=2$, i.e., $l_1$ and~$l_2$ span a
sublattice $\BU(2)$, and $\conj_*$ acts via
\[\label{eq:1}
l_1\mapsto-l_2,\qquad l_2\mapsto-l_1.
\]
Each (simple) singular point of~$C$ gives rise to a singular point
of the double covering, and the exceptional divisors of its
resolution span a root system in $H_2(X;\ZZ)$ of the same type
($\BA$, $\BD$, or~$\BE$) as the original singular point. These
root systems are orthogonal to each other and to~$l_1$, $l_2$;
denote their sum by~$\Sigma$. If all singular points are real,
then $\conj_*$ acts on~$\Sigma$ via multiplication by~$(-1)$.

\begin{lemma} The sublattice $\Sigma\subset H_2(X;\ZZ)$
is quasi-primitive in its primitive hull.


\end{lemma}

\begin{proof} Let $r\notin\Sigma$ be a root in the primitive hull
of~$\Sigma$. Since, obviously, $\Sigma\subset\Pic X$ and
$H_2(X;\ZZ)/\!\Pic X$ is torsion free, one has $r\in\Pic X$. Then,
the Riemann-Roch theorem implies that either~$r$ or~$-r$ is
effective, i.e., it is realized by a $(-2)$-curve in~$X$
(possibly, reducible), which is not contracted by the blow down
(as $r\notin\Sigma$). On the other hand, $r$ is orthogonal
to~$l_1$ and~$l_2$. Hence, the curve projects to a curve in~$Y$
orthogonal to both the rulings, which is impossible.
\end{proof}

\subsection{The calculation}\label{par:2}

\begin{lemma}\label{lem:3.1} The lattice $\Sigma=3 \BA_2$ has no
non-trivial quasi-primitive extensions.
\end{lemma}

\begin{proof} Up to automorphism of~$3\BA_2$, the discriminant
group $\discr 3\BA_2\cong3\left<-\frac23\right>$ has a unique
isotropic element, which is the sum of all three generators. Then,
for the corresponding extension~$\tilde\Sigma\supset\Sigma$ one
has $\discr\tilde\Sigma=\left<\frac23\right>$, i.e.,
$\tilde\Sigma$ has the genus of~$\BE_6$.
Since the latter is unique in its genus (see \cite{Nikulin}), one
has $\tilde\Sigma\cong\BE_6$. Alternatively, one can argue that,
on one hand, an imprimitive extension of~$3\BA_2$ is unique and,
on the other hand, an embedding $3\BA_2\subset\BE_6$ is known: if
$2\BA_2$ is embedded into~$\BE_6$ via the Dynkin diagrams, the
orthogonal complement is again a copy of~$\BA_2$.
\end{proof}

\begin{lemma}\label{lem:3.2} Up to automorphism, the lattice
$\Sigma=8\BA_2$ has two non-trivial quasi-primitive extensions
$\tilde\Sigma\supset\Sigma$; one has $\ell_3(\tilde\Sigma)=6$
or~$4$.
\end{lemma}

\begin{proof} We will show that there are at most two classes.
The fact that the two extensions constructed are indeed
quasi-primitive is rather straightforward, but it is not needed in
the sequel.

Let $\CS=\discr\Sigma\cong8\left<-\frac23\right>$ be the
discriminant group, and let~$G$ be the set of generators of~$\CS$.
The automorphisms of~$\Sigma$ act via transpositions of~$G$ or
reversing some of the generators. (Recall that the decomposition
of a definite lattice into an orthogonal sum of indecomposable
summands is unique up to transposing the summands.) For an element
$a\in\CS$ define its {\it support\/} $\supp a\subset G$ as the
subset consisting of the generators appearing in the expansion
of~$a$ with a non-zero coefficient.
Since each nontrivial summand in the expansion of an element
$a\in\CS$ contributes $-\frac23\bmod2\ZZ$
to the square, $a$~is isotropic if and only if
$\#\supp a=0\bmod3$;
in view of Lemma~\ref{lem:3.1}, such an element cannot belong
to the kernel of a quasi-primitive extension unless $\#\supp a=6$.
(Indeed, if $\#\supp a=3$, then $a$ belongs to the discriminant
group
of
the sum~$\Sigma'$ of certain three of the eight $\BA_2$-summands
of~$\Sigma$, and already~$\Sigma'$ is not primitive, hence, not
quasi-primitive.)

All elements $a\in\CS$ with $\#\supp a=6$ form a single orbit of
the action of $\Aut\Sigma$, thus giving rise to a unique
isomorphism class of quasi-primitive extensions
$\tilde\Sigma\supset\Sigma$ with $\ell_3(\discr\tilde\Sigma)=6$.
Consider the extensions with $\ell_3(\discr\tilde\Sigma)=4$, i.e.,
those whose kernel~$\CK$ is isomorphic to $\ZZ_3\oplus\ZZ_3$. Up to the
action of $\Aut\Sigma$ the generators $g_1,\ldots,g_8$ of~$\CS$
and two elements $a_1$, $a_2$ generating~$\CK$
can be chosen so that $a_1=g_1+\ldots+g_6$ and
$a_2=(g_1+\ldots+g_p-g_{p+1}-\ldots-g_{p+q})+\sigma$, where
$\sigma=0$, $g_7$, or $g_7+g_8$ and $p\ge q\ge0$ are certain
integers
such that
$p+q=\#(\supp a_1\cap\supp a_2)\le6$.
Since $\supp a_1$ and $\supp a_2$ are two six element sets
and $\#(\supp a_1\cup\supp a_2)\le8$, one has
$p+q\ge4$.
Furthermore,
since
$a_1\cdot a_2=\frac23(p-q)\bmod\ZZ=0$,
one has $p-q=0\bmod3$. This leaves three
pairs of values: $(p,q)=(2,2)$, $(3,3)$, or $(4,1)$. In the first
case, $(p,q)=(2,2)$, one does obtain a quasi-primitive extension,
unique up to automorphism. In the other two cases one has
$\#\supp(a_1-a_2)=3$ and, hence, the extension is not
quasi-primitive due to Lemma~\ref{lem:3.1} (cf. the previous
paragraph).

Note that, in the only quasi-primitive case $(p,q)=(2,2)$, for any
pair~$a_1$, $a_2$ of generators of~$\CK$
one has
\begin{equation}\label{eq:2} \supp a_1\cup\supp
a_2=G\quad\text{and}\quad \#(\supp a_1\cap\supp a_2)=4.
\end{equation}

As a
by-product, the same relations must hold for any two independent
(over~$\ZZ_3$) elements~$a_1$, $a_2$ in the kernel of any
quasi-primitive extension.

Now, assume that the kernel of the extension
$\tilde\Sigma\supset\Sigma$ contains
$\ZZ_3\oplus\ZZ_3\oplus\ZZ_3$, i.e.,
$\ell_3(\discr\tilde\Sigma)<4$. Pick three independent
(over~$\ZZ_3$) elements $a_1$, $a_2$, $a_3$ in the kernel. In view
of~\eqref{eq:2}, the principle of inclusion and exclusion implies
that $\#(\supp a_1\cap\supp a_2\cap\supp a_3)=2$. Important is the
fact that the intersection is nonempty. Hence, with appropriate
choice of the signs, there is a generator of~$\CS$, say,~$g_1$,
whose coefficients in the expansions of all three elements~$a_i$
coincide. Then the two differences $b_1=a_1-a_3$ and $b_2=a_2-a_3$
belong to the kernel, are independent, and their supports do not
contain~$g_1$. This contradicts to~\eqref{eq:2}.
\end{proof}

\begin{proposition}\label{prop:3.3} Let $L$ be a lattice isomorphic to
$2\BE_8\oplus3\BU$, and let $S=\Sigma\oplus\BU(2)$ be a sublattice
of~$L$ with $\Sigma\cong8\BA_2$ quasi-primitive in its primitive
hull. Then $L$ has no involutive automorphism~$c$ acting
identically on~$\Sigma$, interchanging the two elements of a
standard basis of~$\BU(2)$, and having exactly two positive
squares in the $(+1)$-eigenlattice $L^{+c}=\ker(1-c)\subset L$.
\end{proposition}

\begin{proof} Assume that such an involution~$c$ exists.
Let~$\tilde\Sigma$ and $\tilde S$ be the primitive hulls
of~$\Sigma$ and~$S$, respectively, in~$L$, and let $T=S^\perp$ be
the orthogonal complement. 
The lattice~$T$ has rank~$4$ and
signature~$0$, {\it i.e.}, it has two positive and two negative squares.

Since $\discr\BU(2)=\ZZ_2\oplus\ZZ_2$ (as a group) has $2$-torsion
only, the $3$-torsion parts $(\discr\tilde\Sigma)\otimes\ZZ_3$ and
$(\discr\tilde S)\otimes\ZZ_3$ coincide. In particular, $c$~must
act 
identically on $(\discr\tilde
S)\otimes\ZZ_3$ (as, by the assumption, so it does on~$\Sigma$)
and, hence,
on $(\discr T)\otimes\ZZ_3$, see
Theorem~\ref{thm:Nik2}.
Furthermore, due to Lemma~\ref{lem:3.2} one has
$\ell_3(\discr\,T)=\ell_3(\discr\tilde S)\ge4$. On the other hand,
$\ell_3(\discr\,T)\le\rank T=4$. Hence, $\ell_3(\discr T)=\rank
T=4$ and the canonical homomorphism
$T^\spcheck\otimes\ZZ_3\to(\discr\,T)\otimes\ZZ_3$ is an
isomorphism. Thus, $c$~must also act identically on
$T^\spcheck\otimes\ZZ_3$ and, hence, both on~$T^\spcheck$ and
$T\subset T^\spcheck$. Indeed, for any free abelian group~$V$, any
involution $c:V\to V$, and any odd prime~$p$, one has a direct sum
decomposition
$V\otimes\ZZ_p=(V^{+c}\otimes\ZZ_p)\oplus(V^{-c}\otimes\ZZ_p)$.
Hence, $c$ acts identically on~$V$ (i.e., $V^{-c}=0$) if and only
if it acts identically on $V\otimes\ZZ_p$ (i.e.,
$V^{-c}\otimes\ZZ_p=0$).

It remains to notice that, under the assumptions, the
skew-invariant part
$S^{-c}=\ker(1+c)\cong\Sigma\oplus\left<-4\right>$ is negative
definite.
Since the total skew-invariant part $L^{-c}$ has exactly one
($=3-2$) positive square, one
of the two positive squares of~$T$,
should fall
to~$T^{-c}$
and the other,
to~$T^{+c}$
In particular, $T^{-c}\ne0$, and the action of~$c$
on~$T$ is not identical.
\end{proof}

Now we have finally reached the goal of this section.

\begin{theorem} The ellipsoid ${\bf Ell}=(Y,s)$
(see~\S\ref{intro} and  \S\ref{par:2}) does not contain a real curve~$C$ of
bi-degree~$(4,4)$ having eight real cusps (and no other
singularities).
\end{theorem}

\begin{proof} Any such curve~$C$ would be irreducible; hence, as
in~\S\ref{par:2}, it would give rise to a sublattice
$8\BA_2\oplus\BU(2)\subset L=H_2(X;\ZZ)\cong2\BE_8\oplus3\BU$ and
involution $c=-\conj_*:L\to L$ which do not exist due to
Proposition~\ref{prop:3.3}.
\end{proof}

\section {Remarks and problems}
\label{sc:rmk}

\noindent I.  Analogously to the  total reality property for
rational curves one can ask a similar question for projective
curves of any genus, namely

\begin{problem}
    Given a real algebraic curve $(\C,\sigma)$ with compact $\C$ and nonempty real part
    $\C_{\sigma}$ and a complex algebraic map $\Psi:\C\to\cP^n$ such that the
inverse images of all the flattening points of $\Psi(\C)$ lie on the real
part  $\C_{\sigma}\subset \C$ is it true that  $\Psi$  is a real algebraic
 up to a projective automorphism  of the image $\cP^n$?
  \label{conj:extSS}
    \end{problem}

    The feeling is that this problem has a negative answer.

\medskip
 \noindent
 II. In the recent \cite{EGSV} the authors found another generalization of
 the conjecture on total reality in case of the usual rational functions.

 \begin{problem}  Extend the results of \cite{EGSV} to  the case of meromorphic
 functions  on curves of
 higher genera.
 \end{problem}

\end{document}